\documentclass{amsart}

\usepackage{mathtools}
\usepackage{amssymb, hyperref}
\usepackage{amsmath,amssymb,amsbsy,amsfonts, latexsym,amsopn,amstext, amsxtra,euscript,amscd}
\usepackage{amsthm}
\usepackage[usenames,dvipsnames]{color}
\usepackage{amssymb, amsthm, amsmath,  amssymb, amsbsy,   amsfonts,  latexsym,  amsopn,   amstext, amsxtra,  euscript,   amscd}
\usepackage{cite}

\usepackage{cite}
\usepackage[shortlabels]{enumitem}

\usepackage[capitalise]{cleveref}

\usepackage{amsrefs}

\hypersetup{
  colorlinks   = true, 
  urlcolor     = blue, 
  linkcolor    = blue, 
  citecolor   = red 
}

\newtheorem{thm}{Theorem}
\newtheorem{prop}{Proposition}
\newtheorem{lem}{Lemma}
\newtheorem{cor}{Corollary}
\newtheorem{rem}{Remark}
\newtheorem{exa}{Example}
\newtheorem{prob}{Problem}

\theoremstyle{definition}
\newtheorem{defi}[thm]{Definition}
\theoremstyle{remark}

\crefname{thm}{Thm.}{}
\crefname{prop}{Prop.}{}
\crefname{lem}{Lem.}{}
\crefname{cor}{Cor.}{}
\crefname{table}{Table}{}

\def\p{\mathfrak p}
\DeclareMathOperator\wgcd{\mathit{wgcd }}   
\DeclareMathOperator\wh{\mathfrak{h}}   
\DeclareMathOperator\awh{\mathfrak{\tilde h} }   

\def\awgcd{\overline{\wgcd }}
\newcommand{\norm}[1]{\left\lVert\mspace{1mu}#1\mspace{1mu}\right\rVert}

\DeclareMathOperator\GL{GL}
\DeclareMathOperator\SL{SL}
\newcommand\w{\mathbf{w}}
\newcommand\x{\mathbf{x}}

\def\0{\mathbf0}

\def\a{{\alpha }}
\def\b{{\beta }}

\begin{document}

\title{On the semistability of binary forms over number fields}
\author{Elira Curri}



 
\begin{abstract}
Let $K$ be a number field, ${\mathcal O}_K$ its ring of integers, and   $ f(x, y) \in {\mathcal O}_K[x, y]$    a binary form with integer coefficeints.
For any given prime $p \in {\mathcal O}_K$ we determine explicitly a binary form $g$   (resp. $\bar f$), $\GL_2 (K)$-equivalent to $f$ which is semistable over the local field $K_p$ (resp. the global field $K$).  Moreover, if  $\xi(f) $ is the corresponding moduli point in the weighted projective space ${\mathbb{WP}}_\w^n (K)$ for a strictly semistable binary form $f$, we determine  the weighted moduli height $\wh (\xi(f))$ for $d=4, 6, 8, 10$.
\end{abstract}

\maketitle


\section{Introduction}
Let $K$ be a number field,  ${\mathcal O}_K$ its ring of integers, and $f\in K [x, y]$ a degree $d\geq 2$ binary form.   
The equivalence classes of binary forms $f (x, y)$, over some algebraic closure of $K$,  are determined by the set of generators of the ring of invariants ${\mathcal R}_d$ of degree $d$ binary forms. It is well known that ${\mathcal R}_d$ is finitely generated.  Let $ \xi_0, \dots ,  \xi_n$ be such generators.   They are homogenous polynomial of degree  $\deg \xi_i = q_i$.   We denote by  $ \xi(f)= \left( \xi_0 (f), \dots ,  \xi_n (f)\right)$ the values of such invariants evaluated at the     form  $f(x, y)$. 
The equivalence class of $f$ is determined by the moduli point  $\xi(f) = [\xi_0 (f): \dots : \xi_n (f)]$  in the weighted projective space ${\mathbb{WP}}_\w^n (K)$ with weights $\w=(q_0, \dots , q_n)$; see   \cite{b-g-sh} for details.    

A   binary form has a root of multiplicity greater than $\frac d 2$  if and only if $\xi (f) = \0$  (cf. \cref{lem-0}).  
Using Hilbert-Mumford numerical criterion (cf. \cref{thm-2}) one proves that   $f(x, y)$ is semistable (resp. stable) if and only if it has no root of multiplicity $\geq \frac d 2$   (resp. $ > \frac d 2$). Hence, for a prime $p \in {\mathcal O}_K$, a binary form is semistable over the completion $K_p$  if and only if $p \nmid \xi_i$, for all $i=0, \dots , n$.  

The main focus of this paper is to find a semistable model for any  given $f \in {\mathcal O}_K[x, y]$. For a fixed prime $p \in {\mathcal O}_K$,   $f(x, y)$ is semistable over the completion $K_p$  if and only if $p \nmid \gcd \left(\xi_0(f), \dots , \xi_n (f) \right)$.   By taking a minimal model of $f(x, y)$ as described in \cite{reduction} we can assume that the weighted greatest common divisor $\wgcd (\xi (f))=1$; see \cite{b-g-sh}.  This is equivalent to assume that $f(x, y)$ can be unstable only for primes $p$ dividing $\frac {\gcd (\xi (f) )} {\wgcd (\xi (f))}$. For each such prime $p \in {\mathcal O}_K$ 
 we find a matrix $A_p = \begin{bmatrix} a & b \\ c & d \end{bmatrix}    \in \GL_2 (K)$ such that   $f^{A_p} (x, y)= f(ax + by, cx + dy) $ is semistable over $K_p$.  
 Taking the matrix  $M=\prod  A_p$ for $p$ dividing $\frac {\gcd (\xi (f) )} {\wgcd (\xi (f))}$ gives a binary form $f^M(x, y) \in {\mathcal O}_K [x, y]$ which is 
    semistable over a global field $K$.

Another goal is to determine  stability in terms of heights in the moduli space ${\mathbb{WP}}_\w^n (K)$.  A natural height in ${\mathbb{WP}}_\w^n (K)$ is the weighted moduli height as defined in \cite{b-g-sh}.  If $f(x, y)$ is unstable than  $\xi (f) = \0$.  When $f(x, y)$ is strictly semistable it has a root of multiplicity $\frac d 2 $.  Such forms do not exist when $d$ is odd and there is only one such form (up to equivalence) for $d$ even.  We compute the weighted moduli height $\wh (\xi (f)$ for degree $d=4, 6, 8, 10$. Moreover we show that $\wh (\xi(f)) \geq 1$ for all semistable forms. 

In \cite{zhang} Zhang defined the \textit{invariant height} (also known as GIT height) and considered the question of determining the semistability (resp. stability) in terms of such height. Such question was considered in more detail for binary forms in \cite{rabina} and \cite{r-savin}, where the authors show that such height is bounded from below for semistable points and give some lower bounds for cubic binary forms.   It is unclear how the invariant height in \cite{zhang} relates to the weighted moduli height in \cite{b-g-sh}.  The comparison between the invariant height and weighted heights seems to raise many questions, which we intend to explore in \cite{c-sh-2}. 
 
This paper is organized as follows. In \cref{sec-2} we give a brief review of the basics in invariant theory and weighted projective spaces. We display generating invariants for the ring of invariants ${\mathcal R}_d$ for $ d \geq 3$ and  $d \leq 10$. 

In \cref{sec-3} are covered some of the classical results of semistability of binary forms including the Hilbert-Mumford criteria. While this material is part of the classical Geometric Invariant Theory, it provides an outline here to prove the fact that a degree $d$ binary form is semistable (resp. stable) if and only if it has a root of multiplicity $\leq \frac d 2$ (resp. $< \frac d 2$). Moreover, we consider binary forms over a number field $K$.  For every prime $p \in {\mathcal O}_K$   we give a condition for a binary form to be semistable over $K_p$  in terms of the coordinates of the point in the weighted moduli point and provide a method how to determine an equivalent form to the given form which is semistable over $K_p$. 

In \cref{sec-4} we define the weighted height and consider the weighted height of   semistable points and strictly semistable points.     
From \cref{semistable-0}  for $d$ odd that are no  strictly semistable binary forms  and for $d$ even there is exactly  one such binary form (up to equivalence).  We display the corresponding points in the weighted moduli space  for 
such binary forms and their corresponding weighted height for $d=4, 6, 8, 10$; see \cref{tab-1}.  

One obvious observation from \cref{thm-6} and  \cref{tab-1} seems that the   weighted height $\wh (\xi(f))$ seems to be growing fast as $d$ increases. This seems to be different from the behavior  of the invariant height and the results in \cite{rabina}.  There is no known estimate for 
$ \displaystyle \lim_{d\to \infty} \frac   {\wh (\xi(f))} d$.

It would also be interesting to estimate the number of stable binary forms with integer coefficients such that their weighted moduli height is less than the weighted moduli height of the strictly semistable binary form of the same  even degree $d$. Equivalently this would estimate the number of binary forms (up to $K$-equivalence) with weighted height less than the strictly semistable form, such that the field of moduli is also a field of definition.  Some of these problems are intended to be discussed in \cite{c-sh-2}. 

\section{Preliminaries}\label{sec-2}
Let  $k$ be a field,  $k[x, y]$  be the  polynomial ring in  two variables and   $V_d$ denote  the $(d+1)$-dimensional  subspace  of  $k[x, y]$  consisting of homogeneous polynomials
\begin{equation}\label{eq1}
f(x,y) = a_d x^d + a_{d-1} x^{d-1}y +  \dots  + a_0 y^d
\end{equation}
of  degree $d$. Elements  in $V_d$  are called  \textbf{binary  forms} of degree $d$.    $\GL_2(k)$ acts as a natural group of automorphisms on $k[x, y] $.     Denote by $f \to f^M$ this action.  It is well  known that $\SL_2(k)$ leaves a bilinear  form (unique up to scalar multiples) on $V_d$ invariant; see \cite{reduction} for details.  

If $k$ is algebraically closed then  $f(x,y)$  can be factored as
\begin{equation} \label{eq4}
f(x,y)  = (\b_1  x  -  \a_1 y) \cdots (\b_d  x  - \a_d  y) = \displaystyle \prod_{1 \leq  i \leq  d} \det
\left(\begin{pmatrix} x & \a_{i} \\ y & \b_i
\end{pmatrix} \right)
\end{equation}
The points  with homogeneous coordinates $(\a_i, \b_i)  \in \mathbb P^1$ are  called the  \textbf{projective roots}  of  $f$.
For $M  \in \GL_2(k)$ we have
\[
 f^M (x,y)   = \left( \det  M \right)^{d}  (\b_1^{'}  x -  \a_1^{'}  y) \cdots (\b_d^{'} x  - \a_d^{'} y),
\]
where
$\begin{pmatrix}   \a_i^{'}  \\  \b_i^{'}  \end{pmatrix} = M^{-1}  \begin{pmatrix} \a_i\\ \b_i \end{pmatrix}$.

Consider $a_0, a_1,  \ldots, a_d$ as transcendentals over $k$  (coordinate  functions on $V_d$). Then the coordinate  ring of $V_d$ can be identified with $ k[a_0  , \ldots ,  a_d] $. We define an action of $\GL_2(k)$ on $k[a_0, \dots , a_d]$ via
\[
\begin{split}
\GL_2 (k) & \times k[a_0, \ldots, a_d] \to  k[a_0, \ldots, a_d]  \\
(M, F)   & \to F^M : = F (f^M), \quad \text{ for all } \; f \in V_d.
\end{split}
\]
Thus for a polynomial $F \in k[a_0, \ldots ,  a_d]$ and $M \in \GL_2(k)$, define $F^M \in k[a_0, \dots , a_d]$ as 
${F^M}(f):= F ( f^M),$
for all $f \in V_d$. Then  $F^{MN} = (F^{M})^{N}$. 
The homogeneous degree in $a_0, \dots , a_d$ is called the \textbf{ degree} of $F$,  and the homogeneous degree in $x, y$ is called the \textbf{  order} of $F$.   An \textbf{invariant}  is usually referred to  an $\SL_2(k)$-invariant on $V_d$.
 Hilbert's theorem  says that the ring of invariants ${\mathcal R}_d$ of binary forms of degree $d$ is finitely generated.    Thus, 
 \begin{enumerate}[\upshape(i)]
\item ${\mathcal R}_d$ is finitely generated, and 
\item ${\mathcal R}_d$ is a graded ring.
\end{enumerate}

Let $ \xi_0, \dots  ,  \xi_n$ be a minimal set of  generators of ${\mathcal R}_d$.  Since $\xi_i\in k[a_0, \dots , a_d]$ are homogenous polynomials we denote 
$\deg \xi_i= q_i$.   
The set of degrees $(q_0, \dots ,q_n)$ is often called the \textbf{set of weights}.

The following result is crucial for the rest of this paper   (see  \cite[Prop.~1]{reduction}  for its proof). 

\begin{lem}\label{lem-1}
Let $f,g\in V_d, M\in \GL_2 (k),  \l=\left(\det M\right)^{\frac d 2}$. Then  $f=g^M$     if and only if 
\[ 
\left(  \xi_0 (f), \dots    \xi_i (f), \dots ,  \xi_n (f) \right) = \left( \l^{q_0} \,  \xi_0 (g), \dots ,    \l^{q_i}\,   \xi_i (g), \dots , \l^{q_n} \,  \xi_n (g)    \right).
\] 
\end{lem}

While the proof of the above result is technical and requires some knowledge of classical theory of invariants, its understanding it is rather elementary.  

\begin{exa}
Take   a binary quadratic
$ f(x, y) = a_2 x^2 + a_1 xy + a_0 y^2$.
${\mathcal R}_2$ is generated by the discriminant ${\Delta}=a_1^2 - 4a_0a_2$, which is an   $SL_2(k)$-invariant.  For any matrix 
$M=\begin{bmatrix} a & b \\ c & d \end{bmatrix} \in \GL_2 (k)$  the form $f^M (x, y):= f(ax+by, cx+dy) $
has discriminant  
\[
{\Delta} (f^M) = (ad-bc)^2 \cdot (a_1^2 - 4a_0a_2)=   \left( \det M \right)^2 \cdot {\Delta}(f).
\] 
In general, the exponent of $\det (M)$ is precisely $\frac d 2 \cdot q_i$, where $q_i$  the degree of the invariant $\xi_i$.  
\end{exa}

\begin{lem}\label{lem-0}    
If $k={\mathbb Q}$  we can choose  $ \xi_0, \dots  ,  \xi_n$ with integer coefficients and   primitive polynomials in $\mathbb Z [a_0, \dots , a_d]$
(i.e. the greatest common divisor of coefficients of each $\xi$ is 1).  
\end{lem}

\proof
Without loss of generality we can assume that each $f \in {\mathbb Q} [x, y]$ has integer coefficients since binary forms are defined up to multiplication by a constant.   Then for $\xi = [ \xi_0, \ldots , \xi_n]$ each coordinate   $\xi_i (f) \in {\mathbb Q}[a_0, \ldots , a_d]$.   Then multiplying the tuple 
$\xi = \left[  \xi_0 : \cdots : \xi_n     \right] $
by the least common multiple of all the denominators, say $\lambda$, we get  a representative 
\[
\left[  \l^{q_o} \xi_0 : \cdots : \l^{q_i} \xi_i, \cdots , \l^{q_n}  \xi_n     \right] 
\]
of $\xi $ with integer coefficients.  We can redefine each $\xi_i$ by taking its primitive part. 

\qed


\subsection{$\mbox{Proj } {\mathcal R}_d$ as a weighted projective space} 
Let $ \xi_0, \dots  ,  \xi_n$ be the generators of ${\mathcal R}_d$ with degrees $q_0, \dots , q_n$ respectively.  
Since all $ \xi_0, \dots ,  \xi_i, \dots ,  \xi_n$ are homogenous polynomials then ${\mathcal R}_d$ is a graded ring and  $\mbox{Proj } {\mathcal R}_d$ as a weighted projective space. 
Let   $\w:=(q_0, \dots , q_n) \in \mathbb Z^{n+1}$ be a fixed tuple of positive integers called \textbf{weights}.   Consider the action of $k^\star = k \setminus \{0\}$ on ${\mathbb A}^{n+1} (k)$ as follows
\[ \lambda \star (x_0, \dots , x_n) = \left( \l^{q_0} x_0, \dots , \l^{q_n} x_n   \right) \]
for $\l\in k^\ast$.  The quotient of this action is called a \textbf{weighted projective space} and denoted by   ${\mathbb{WP}}^n_{(q_0, \dots , q_n)} (k)$. 
It is the projective variety $Proj \left( k [x_0,...,x_n] \right)$ associated to the graded ring $k [x_0, \dots ,x_n]$ where the variable $x_i$ has degree $q_i$ for $i=0, \dots , n$. 
We denote greatest common divisor of $q_0, \dots , q_n$ by $\gcd (q_0, \dots , q_n)$.    The space ${\mathbb{WP}}_w^n$ is called \textbf{well-formed} if   
\[ \gcd (q_0, \dots , \hat q_i, \dots , q_n)   = 1, \quad \text{for each } \;  i=0, \dots , n. \]
While most of the papers on weighted projective spaces are on well-formed spaces, we do not assume that here.   
We will denote a point $\p \in {\mathbb{WP}}_w^n (k)$ by $\p = [ x_0 : x_1 : \dots : x_n]$. 

\begin{lem} Let $ \xi_0,   \xi_1,   \dots,   \xi_n$ be the generators of the ring of invariants ${\mathcal R}_d$ of degree $d$ binary forms. 
A $k$-isomorphism class of a binary form $f$ is determined by  the point 
\[  \xi (f) := \left[  \xi_0 (f),  \xi_1 (f), \dots ,  \xi_n(f) \right]  \in {\mathbb{WP}}_\w^n  (k). \]
Moreover, for any two forms $f,$ and $g$ we have that 
$f=g^M$ for some $M\in \GL_2 (k)$ if and only if  $  \xi(f)  = \lambda\star  \xi(g)$,    for $\lambda = \left( \det A \right)^{\frac d 2}$.
\end{lem}

\proof
The proof is a direct consequence of \cref{lem-1}. 
\qed

\subsection{Generating invariants}
Finding generators for the ring of invariants $R_d$ is a classical problem in which worked many of the most important mathematicians of the XIX-century.  Such invariants are generated in terms of transvections or root differences.  
%
For given binary invariants $f, g \in V_d$  the $r$-th transvection of $f$ and $g$, denoted by $(f, g)_r$,  is  defined as  
%
\[(f,g)_r:= \frac {(m-r)! \, (n-r)!} {n! \, m!} \, \,
\sum_{k=0}^r (-1)^k
\begin{pmatrix} r \\ k
\end{pmatrix} \cdot
\frac {\partial^r f} {\partial x^{r-k} \, \,  \partial y^k} \cdot \frac {\partial^r g} {\partial x^k  \, \, \partial y^{r-k} },
\]
see Grace and Young \cite{GY} for  details.    
Transvections are a convenient way of generating invariants since they are expressed in terms of the coefficients of the binary form, on contrary to the method of generating invariants through root differences which give the invariants in terms of roots of binary forms, which we will briefly describe next.

While there is no method known to determine a generating set of invariants for any ${\mathcal R}_d$,  we display a minimal generating set for all $3\leq d \leq 10 $.  For the rest of this section  $f(x, y)$ is given as in \cref{eq1} and a minimal set of invariants is always picked as in \cref{lem-0}.
\subsubsection{Cubics}     A generating set for ${\mathcal R}_3$ is   $\xi = \{ \xi_0 \}$, where 
 \[
 \xi_0= \left(  (f, f)_2, (f, f)_2    \right)_2
 \]
  and 
\[
\xi_0 = -54 a_0^2 a_3^2+36 a_1 a_3 a_0 a_2-8 a_2^3 a_0-8 a_1^3 a_3+2 a_2^2 a_1^2
\]
\subsubsection{Quartics}   A generating set for ${\mathcal R}_4$ is   $\xi = [ \xi_0 , \xi_1]$ with $\w=(2,3)$, where 
%
\[
\xi_0 =   (f, f)_4   \; \text{ and } \; \xi_1  = \left(   f, (f, f)_2   \right)_4
\]
In terms of $a_0, \dots , a_4$, the invariants   $\xi_0$ and $\xi_1$  are 
\[
\begin{split}
\xi_0 & =12 a_0 a_4-3 a_1 a_3+a_2^2 \\
\xi_1 & =72 a_0 a_2 a_4-27 a_0 a_3^2-27 a_1^2 a_4+9 a_1 a_2 a_3-2 a_2^3 \\
\end{split}
\]

\subsubsection{Quintics}   A generating set for ${\mathcal R}_4$ is   $\xi = [ \xi_0 , \xi_1, \xi_2]$ with $\w=(4, 8, 12)$, where  
\[
\xi_0=   (c_1, c_1)_2, \quad \xi_1  = (c_4, c_1)_2, \quad \xi_2 = (c_4, c_4)_2, 
\]
for   
\[
c_1 = (f, f)_ 4, \; c_2 = (f, f)_2, \; c_3= (f, c_1)_2, \; c_4 = (c_3, c_3)_2.
\] 
We don't display such invariants in terms of $a_0, \dots , a_5$ since their expressions are big.

\subsubsection{Sextics}   The case of sextics was studied in detail due to their connection to genus 2 curves.  Generating sets were known in detail
by XIX-century mathematicians   (Bolza, Clebsch, et al.) when $\mbox{char } k =0$ and by Igusa for $\mbox{char } k >0$.  For a more modern treatment see \cite{k-sh-v} or    \cite{m-sh-1, m-sh-2, m-sh-3}, 
where    invariants of binary sextics are used to study the moduli space of genus 2 curves and  even expressed  in terms of modular forms. 
To have a uniform treatment of invariants in this paper we will define a generating set for the case $d=6$ slightly different from generating sets commonly used in literature.  

Let $c_1 = (f, f)_ 4$, $c_3 = (f, c_1)_4$, $c_4 = (c_1, c_1)_2$.  A generating set for ${\mathcal R}_6$ is   $\xi = [ \xi_0 , \xi_1, \xi_2, \xi_3]$ with weights $\w=(2,4,6, 10)$   (we are assuming $\mbox{char } k \neq 2$), where 
\[
\xi_0=   (f, f)_6, \; \xi_1= (c_1, c_1)_4, \;   \xi_2   = (c_4, c_1)_4, \; \xi_3  = (c_4, c_3^2)_4 
\] 
and
\begin{small}
\[
\begin{split}
\xi_0 & = 120 a_0 a_6-20 a_1 a_5+8 a_2 a_4-3 a_3^2 \\
\xi_1 & = 7500 a_0^2 a_6^2-2500 a_1 a_6 a_0 a_5-200 a_2 a_6 a_0 a_4+500 a_2 a_5^2 a_0\\
& +300 a_0 a_6 a_3^2-300 a_3 a_5 a_0 a_4+80 a_4^3 a_0+500 a_1^2 a_6 a_4-300 a_2 a_6 a_1 a_3 \\
& -100 a_2 a_5 a_1 a_4+100 a_3^2 a_5 a_1-20 a_4^2 a_1 a_3+80 a_2^3 a_6-20 a_3 a_5 a_2^2+28 a_4^2 a_2^2 \\
& -16 a_3^2 a_4 a_2+3 a_3^4 \\
\xi_2 & = -125000 a_0^3 a_6^3+62500 a_0^2 a_1 a_5 a_6^2+35000 a_0^2 a_2 a_4 a_6^2-25000 a_0^2 a_2 a_5^2 a_6-7500 a_0^2 a_3^2 a_6^2\\
& -7500 a_0^2 a_3 a_4 a_5 a_6+6250 a_0^2 a_3 a_5^3+4000 a_0^2 a_4^3 a_6-2500 a_0^2 a_4^2 a_5^2-25000 a_0 a_1^2 a_4 a_6^2\\
& -7500 a_0 a_1 a_2 a_3 a_6^2 +10000 a_0 a_1 a_2 a_4 a_5 a_6+6250 a_0 a_1 a_3^2 a_5 a_6-3500 a_0 a_1 a_3 a_4^2 a_6\\
& -2500 a_0 a_1 a_3 a_4 a_5^2+1000 a_0 a_1 a_4^3 a_5+4000 a_0 a_2^3 a_6^2 -3500 a_0 a_2^2 a_3 a_5 a_6-600 a_0 a_2^2 a_4^2 a_6\\
& +1100 a_0 a_2 a_3^2 a_4 a_6+250 a_0 a_2 a_3^2 a_5^2+300 a_0 a_2 a_3 a_4^2 a_5-160 a_0 a_2 a_4^4 -150 a_0 a_3^4 a_6+250 a_1^2 a_3^2 a_4 a_6\\
& -150 a_0 a_3^3 a_4 a_5 + 60 a_0 a_3^2 a_4^3+6250 a_1^3 a_3 a_6^2-2500 a_1^2 a_2^2 a_6^2-2500 a_1^2 a_2 a_3 a_5 a_6+250 a_1^2 a_3 a_4^2 a_5\\
& -100 a_1^2 a_4^4+1000 a_1 a_2^3 a_5 a_6+300 a_1 a_2^2 a_3 a_4 a_6+250 a_1 a_2^2 a_3 a_5^2 -100 a_1 a_2^2 a_4^2 a_5-150 a_1 a_2 a_3^3 a_6\\
&-350 a_1 a_2 a_3^2 a_4 a_5+140 a_1 a_2 a_3 a_4^3+100 a_1 a_3^4 a_5-40 a_1 a_3^3 a_4^2-160 a_2^4 a_4 a_6-100 a_2^4 a_5^2\\
& +60 a_2^3 a_3^2 a_6+140 a_2^3 a_3 a_4 a_5-24 a_2^3 a_4^3-40 a_2^2 a_3^3 a_5-8 a_2^2 a_3^2 a_4^2+8 a_2 a_3^4 a_4-a_3^6\\
\end{split}
\]
\end{small}

\begin{rem}
The reader should be aware that usually the invariants of binary sextics are denoted by $[J_2, J_4, J_6, J_{10}]$ with $J_{10}$ being the discriminant of the sextic, but that is not the case here. 
\end{rem}

\subsubsection{Septics}    A generating set of ${\mathcal R}_7$ is given by $\xi= [\xi_0, \xi_1, \xi_2, \xi_3, \xi_4]$ with weights $\w=(4, 8, 12, 12, 20)$.  We define them as follows.   Let
\[
c_1 = (f, f)_ 6, \quad c_2 = (f, f)_ 4, \quad c_4= (f, c_1)_2, \quad  c_5 = (c_2, c_2)_4, \quad c_7 = (c_4, c_4)_ 4   
\]
and 
\[ 
\begin{split}
\xi_0 = (c_1, c_1)_2, \quad  \xi_1 = (c_7, c_1)_2, \quad  \xi_2 = ((c_5, c_5)_2, c_5)_4, \\
 \xi_3  = \left(   (c_4, c_4)_2, c_1^3 \right)_6, \quad  \xi_4 = \left(   \left[(c_2, c_5)_4 \right]^2, (c_5, c_5)_2 \right)_4  
\end{split}
\]
\subsubsection{Octavics}     Finding a generating set for the ring of invariants of binary octavics was one of the biggest achievements of classical invariant theory.  A basis was determined by von Gall; see \cite{vg-1, vg-2}.  Later the ring ${\mathcal R}_8$ was studied by \cite{shioda}, \cite{sh-14} and many others. 
A generating set of ${\mathcal R}_8$ is given by $\xi= [\xi_0, \xi_1, \xi_2, \xi_3, \xi_4, \xi_5]$ with weights $\w=(2,3,4,5,6,7)$.  We define them as follows.  Let 
\[
c_1 = (f, f)_ 6, \quad c_2 = (f, c_1)_ 4, \quad c_3= (f, f)_4, \quad c_5 = (c_1, c_1)_2.
\] 
Then the invariants are:
\[ 
\begin{split}
\xi_0 = (f, f)_8, \quad  \xi_1 = (f, c_3)_8, \quad \xi_2 = (c_1, c_1)_4, \\
 \xi_3 = (c_1, c_2)_4, \quad \xi_5 = (c_5, c_1)_4, \quad \xi_5 = ((c_1, c_2)_2, c_1)_4 
\end{split}
 \]
and 
\[
\begin{split}
\xi_0 & = 280 a_0 a_{8}-35 a_1 a_{7}+10 a_2 a_6-5 a_3 a_5+2 a_4^2 \\
\xi_1 & = 3920 a_{8} a_0 a_4-2450 a_{7} a_0 a_5  +1050 a_0 a_6^2-2450 a_{8} a_1 a_3+735 a_{7} a_1 a_4-175 a_6 a_1 a_5 \\
& +1050 a_{8} a_2^2-175 a_{7} a_2 a_3-110 a_6 a_2 a_4+75 a_2 a_5^2+75 a_6 a_3^2-45 a_5 a_3 a_4+12 a_4^3\\
\xi_3 & = 2458624 a_0^2 a_{8}^2-614656 a_1 a_{8} a_0 a_{7}-12544 a_2 a_{8} a_0 a_6+82320 a_2 a_{7}^2 a_0 \\
& +53312 a_0 a_{8} a_3 a_5-35280 a_3 a_{7} a_0 a_6-25088 a_0 a_{8} a_4^2+4704 a_4 a_5 a_0 a_{7}+8064 a_4 a_6^2 a_0 \\
& -3360 a_5^2 a_0 a_6+82320 a_1^2 a_{8} a_6+2401 a_1^2 a_{7}^2-35280 a_2 a_{8} a_1 a_5-13132 a_2 a_{7} a_1 a_6 \\
& +4704 a_1 a_{8} a_3 a_4+3626 a_3 a_{7} a_1 a_5+3780 a_3 a_6^2 a_1+784 a_1 a_{7} a_4^2-3864 a_4 a_6 a_1 a_5 \\
& +1260 a_5^3 a_1+8064 a_2^2 a_{8} a_4+3780 a_2^2 a_{7} a_5+256 a_2^2 a_6^2-3360 a_2 a_{8} a_3^2-3864 a_3 a_{7} a_2 a_{4 } \\
& -1516 a_3 a_6 a_2 a_5+1984 a_4^2 a_6 a_2-504 a_5^2 a_2 a_4+1260 a_3^3 a_{7}-504 a_4 a_6 a_3^2+589 a_5^2 a_3^{2 } \\
& -320 a_4^2 a_5 a_3+64 a_4^4 \\
\end{split}
\]
\subsubsection{Nonics}       A generating set of ${\mathcal R}_9$ is given by $\xi= [\xi_0, \xi_1, \xi_2, \xi_3, \xi_4, \xi_5, \xi_6]$ with weights $\w=(4, 8, 10, 12, 12, 14, 16)$.   Let 
\[
\begin{split}
& c_1 = (f, f)_ 8, \; c_2 = (f, f)_ 6, \; c_4= (f, f)_2, \; c_5 = (f, c_1)_2, \; c_6 = (f, c_2)_ 6, \\
&  c_7 = (c_2, c_2)_ 4, \; c_9= (c_5, c_5)_4, \; c_{21} = (f, c_2)_2, \; c_{25} = (c_4, c_4)_ {10}, \; c_{27} = (c_6^3, c_6)_ 3 
\end{split}
\]
and 
\[ 
\begin{split}
\xi_0  &= (c_1, c_1)_2, \\  
\xi_1 & = (c_2, c_6^2)_6, \\ 
\xi_2 & = \left( ((c_{25},f)_6, c_{21})_5,c_2 \right)_6, \\  
\xi_3 & =\left ((c_7, c_7)_2,c_7\right)_4 \\
\xi_4  & = (c_9, c_1^3)_6, \\ 
\xi_5 &= \left((c_2, c_{27})_3\right)_6, \\ 
\xi_6 & = \left((c_5, c_5)_2, c_1^5 \right)_{10}.
\end{split}
 \]

\subsubsection{Decimics}     A generating set of ${\mathcal R}_8$ is given by $\xi= [\xi_0, \xi_1, \xi_2, \xi_3, \xi_4, \xi_5, \xi_6, \xi_7, \xi_8]$ with weights $\w=(2, 4, 6, 6, 8, 9, 10, 14, 14)$. 
Let 
%
\begin{align*}
& c_1 = (f, f)_ 8,  		&  c_2 &= (f, f)_ 6,            & 		& c_5 = (f, c_1)_ 4,  	&	&  c_6= (f, c_2)_8, \\
&  c_7 = (c_2, c_2)_6, 	&  c_8 &= (c_5, c_5)_ 4,  &	&  c_9 = (c_2, c_7)_ 4, 	& 	& c_{10}= (c_1, c_1)_2, \\
& c_{16} = (c_5, c_5)_2, 	&  c_{19} &= (c_5, c_1)_ 1, &   & c_{25} = (c_7, c_7)_ 2	&   &  \\            
\end{align*}
%
and 
\[ 
\begin{split}
\xi_ 0 & =  (f, f)_{10}, \\ 
\xi_1 & =  (c_1, c_1)_4, \\
\xi_2 & =  (c_5, c_5)_6, \\ 
\xi_3 & = (c_6, c_6)_2, \\
\xi_4  & =  (c_1, c_8)_4, \\
\xi_5  & =  (c_{19}, c_1^2)_8), \\
\xi_6  & =  (c_{16}, c_1^2)_8, \\
 \xi_7         & =  (c_{25}, c_9)_4,    \\
\xi_8 & =  (c_{10}^2, c_{16})_8, 
\end{split}
 \]
and 
\[
\begin{split}
\xi_1 & = 2520 a_0 a_{10}-252 a_1 a_{9}+56 a_2 a_{8}-21 a_3 a_{7}+12 a_4 a_6-5 a_5^2 \\
\end{split}
\]

\subsection{Root differences}  
Invariants can also be expressed in terms of root differences.  For example the discriminant ${\Delta} (f)= \prod_{i\neq j} (\a_i-\a_j)$. We will not list such formulas here, but the interested reader can check \cite{vishy} among other sources. 

\begin{lem} Let $f\in V_d$.   
\begin{enumerate}[\upshape(i), nolistsep]
\item If $f$ a root of multiplicity $r> \frac d 2$ then $\xi (f) =(\xi_0, \dots , \xi_n)= (0, \dots , 0)$. 
\item If  $d$ is even, then    all binary forms with a root of multiplicity    $\frac d 2$  have the same invariants.  
\end{enumerate}
\end{lem}

\proof   Let $m=\lfloor \frac d 2 \rfloor$.  Hence, $d=2m$ when $d$ is even and $d=2m+1$ when $d$ is odd.   Notice that 
every  degree $s$  invariant $J_s\in k[a_0, \dots , a_d] $   is invariant under the permutation $(a_i, a_{d - i}) $   for $i=0, \cdots , m$, since such permutation corresponds  to permuting $x$ and $y$.  If $f(x, y)$ has a root of multiplicity $r$ then we can assume that $f(x, y)= x^r g(x, y)$  for some degree $(d-r)$ binary form $g(x, y)$, say 
\[  
g(x, y) = b_{d-r} x^{d-r}  + b_{d-r-1} x^{d-r-1} y + \cdots b_1 x y^{d-r-1} + b_0 y^{d-r}.
\]
Then 
\begin{equation}\label{ss-bin}
f(x, y) =  b_{d-r} x^{d}  + b_{d-r-1} x^{d-1} y + \cdots b_1 x^{r+1} y^{d-r-1} + b_0x^{r}  y^{d-r}.
\end{equation}
Hence every $\xi_i (f)$ will be written in terms of coefficients $b_0, \dots , b_{d-r}$ or equivalently in terms of $a_i$, where 
\[
a_0 = \cdots = a_{r-1}=0 \quad \text{ and } \quad a_{r+j} = b_j, \; \text{ for } \; j=0, \dots , d-r.
\]
 Hence, evaluated at $f(x, y)$ as in \cref{ss-bin} is given by  a sum of degree $s$ monomials in $a_r, \dots , a_d$
 since $a_0 = \cdots = a_{r-1}=0$.

To prove part i) let $r= m+1$. Then  all  
$a_0 = \cdots = a_{(m+1)-1}=0 $ which makes   all $a_i=0$, for $i=0, \dots, m $. Hence, $J_s=0$ which implies that 
$(\xi_0 (f), \ldots, \xi_n (f) )=(0, \ldots , 0)$.

To prove ii) let $r=m$ and $d=2m$. Then, all $a_0= \cdots = a_{m-1}=0$. Hence $a_m = b_0$ is the only nonzero coefficient. 
Therefore each $\xi_i (f)$ is a degree $q_i$ homogenous polynomial in $b_0$. Thus,
\[ 
\xi (f)  = \left[     b_0^{q_0} \cdot \lambda_0 : b_0^{q_1} \cdot \lambda_1 : \ldots : b_0^{q_n} \cdot  \lambda_n  \right]
= \left[  \lambda_1 : \cdots : \lambda_n   \right]
\]
for some  $\lambda_i \in k$. Hence,   there is a unique set of invariants 
$\xi (f)  = \left[    \lambda_0 , \lambda_1, \ldots , \lambda_n  \right].$
This completes the proof. 
\qed

%
 
\section{Stability of binary forms}\label{sec-3}

 Let ${\mathcal X}\subset {\mathbb P}^d (k)$, 
 $G \leq \GL_2 (k)$,  and $I  \in k[a_0, \ldots, a_d]$ a $G$-invariant polynomial.   
 By ${\mathcal X}_I \subset {\mathbb P}^d (k) $ we denote the set 
 \[{\mathcal X}_I := \{ \b \in {\mathcal X} \, | \, I(\b) \neq 0\}.\] 
 
\begin{defi} 
A point $\a  \in {\mathcal X}$ is called \textbf{stable under the $G$-action} if $\a$ has a finite stabilizer $G_\a$ and there exist a $G$-invariant  $I\in k[a_0, \ldots, a_d]$ such that  $\a\in {\mathcal X}_I$.   If we drop the condition that the stabilizer $G_\a$ is finite then $\a\in {\mathcal X}$ is called \textbf{semistable under the $G$-action}.
\end{defi} 
 
For any $\a \in  {\mathbb P}^d (k)$ we will denote by $\hat \a$ one of its preimages under the natural projection $\pi: {\mathbb A}^{d+1} (k)  \to {\mathbb P}^d (k)$. 
Let ${\mathbb B}:=\{ b_0, \ldots , b_d \}$ be a basis for ${\mathbb A}^{d+1} (k)=k^{d+1}$ and $\hat \a  \in k^{d+1}$ given by  
\[ 
\hat \a  = \sum_{i=0}^d \hat \a_i b_i.
\]
%

A subgroup $G$ of $\SL_n (k)$ is called a \textbf{1-parameter group} if there is a non-trivial homomorphism of algebraic groups $\l : k^* \to G$.

For $t\in k^*$ let $t. \a$ denote the linear action of $G=k^*$ on ${\mathcal X}$.  
Hence, since $G$ is a one parameter subgroup of $k$ then 
\[
t.\hat \a = \sum_{i=0}^n t^{r_i} \hat \a_i b_i,
\]
see Newstead \cite[pg. 7]{newstead}.     We define 
\[
\mu (\a) := \max \{    - r_i \, | \,  \a_i \neq 0\}  \quad \text{ and } \quad  \mu^- (\a) := \max \{     r_i \, | \,  \a_i \neq 0\} 
\] 
For any  1-parameter subgroup $G$ and a homomorphism  $\l$, sometimes  we write  $\mu (\a, \l )$ for the value of  $\mu (\a)$ for the action of $k^*$ on ${\mathcal X}$ induced by $\l$.

\begin{lem}
For every $\a\in {\mathcal X}$, $\mu(\a)$ is the unique integer such that $\lim_{t\to 0} t^{\mu(\a)} (t.\a)$ exists and is nonzero.  Moreover, $\mu(\a)$ is independent of $\a$ or the basis ${\mathbb B}$  and 
$\mu (\a) > 0$ if and only if $\lim_{t\to 0} t^{\mu(\a)} (t.\a)$ does not exist.
\end{lem}

For a proof see \cite{mumford} or \cite[pg. 7]{newstead} among other places. 

\begin{prop}  
The following hold:
\begin{enumerate}[\upshape(i)]
\item $\a$ is stable if and only if $\mu (\a) >0$ and $\mu^- (\a) > 0$. 
\item  $\a$ is semistable if and only if $\mu (\a) \geq 0$ and $\mu^- (\a)\geq 0$. 
\end{enumerate}
\end{prop}

 
The following result is often used as the definition of stability of binary forms. Its proof is done usually via the Hilbert-Mumford criteria.  

\begin{thm}[Hilbert-Mumford Criterion] \label{thm-2}
The following hold:
\begin{enumerate}[\upshape(i)]
\item  $\a$ is stable if and only if $\mu (\a, \l) >0$ for every 1-PS $\l$ of $G$.
\item  $\a$ is semistable if and only if $\mu (\a, \l) \geq 0$ for every 1-PS   $\l$ of $G$.
\end{enumerate}
\end{thm}

\begin{rem}
For the semistable case when  $G$ is $\SL_n ({\mathbb C})$  there is a proof given by Hilbert using convergent power series. Mumford and Seshadri prooved it for all $k$ and all reductive $G$ using formal power series and a theorem of Iwahori.
\end{rem}
 
 
\begin{thm}
A binary form $f(x, y) \in k [x, y]$ of degree $\deg f=d$  is  stable if and only if all roots  of $f$ are  of multiplicity $< \frac d 2$ 
and semistable if and only if all roots are  of multiplicity $\leq  \frac d 2$.  
\end{thm}

\proof
Any one parameter subgroup $G$ of  $\SL_2 (k)$ is  given by 
\[
\l (t)  = \left\{  \begin{pmatrix}  t^r & 0 \\ 0 & t^{-r} \end{pmatrix}  \, : t\in k^\star\right\}
\]
for some $r\geq 0$.  
Then for $f(x, y)= \sum_{i=0}^d a_i x^i y^{d-i}$ we have 
\[
\l (t) \cdot f(x, y)= \sum t^{r(2i - d)} a_i x^i y^{d-i}.
\]
Hence
\[
\mu (f, \l) = - \min \{  2i - d : a_i \neq 0\} = \max \{ 2i - d : a_i \neq 0\} = 2 i_0 - d,
\]
where $i_0$ is the largest integer for which $\a_i\neq 0$.  Hence, when $\mu (f, \l)\geq 0$
we have $[0, 1]$ as a root  with multiplicity at most $\frac  d 2$ and when 
$\mu (f, \l)> 0$ then $[0, 1]$ has multiplicity strictly less than $\frac d 2$.
This completes the proof.
\qed

\begin{rem}
If   $f(x, y)$ has roots  of multiplicity  $\frac d 2$ we say that $f$ is \textbf{strictly semistable}. 
\end{rem}

\begin{cor}  A binary form $f(x, y)$  of degree $\deg f=d$ 
is unstable if and only if    $ \xi (f)=\0$ in ${\mathbb{WP}}_\w^n (k)$.   
Moreover, if $d$ is even there is only one strictly semistable point in the moduli space and there are no such points when $d$ is odd. 
\end{cor}

\subsection{Stability of binary forms over number fields}\label{sec-5}
Next we want to focus on the stability of binary forms over number fields.  We fix the following notation for the remainder  of this paper. \\


$K$ a number field,

${\mathcal O}_K$ the ring of integers of $K$,

$\nu$ an absolute value of $K$,

$M_K$ the set of all absolute values  of $K$,

$M_{K}^0$ the set of non-Archimedean absolute values of $K$,  

$M_{K}^\infty$ the set of Archimedean absolute values of $K$,

$K_\nu$ the completion of $K$ at $\nu$,

$n_\nu$ the local degree $[K_\nu :{\mathbb Q}_\nu]$ \\


Take $\p =\xi (f) \in {\mathbb{WP}}_\w^n (K)$.   We can assume that $\xi (f)= [\xi_0: \dots : \xi_n]$  has coordinates in ${\mathcal O}_K$. Further assume that  $\p$ is normalized (i.e. $\wgcd (\xi_0, \dots , \xi_n) =1$).  
Let $p$ be a prime in $ {\mathcal O}_K$ such that $p \mid \gcd (\xi_0, \dots , \xi_n)$.   Then $f$ is unstable over $K_p$.  

\subsubsection{Stability over local fields}
 
 Next we show how  to determine an equivalent  binary form $g(x, y)$   to $f(x, y)$  which is semistable over $K_p$. 
 
\begin{lem}\label{semistable-0}
Let $ f \in {\mathcal O}_K [x, y]$ and      $\p=\xi(f)=[\xi_0, \dots , \xi_n]  \in {\mathbb{WP}}_\w^n ({\mathcal O}_K)$.  

\begin{enumerate}[\upshape(i), nolistsep]
\item  $f$  is a semistable binary form over $K_p$      if    and only if        $p \nmid \gcd(\xi_0, \dots, \xi_n)$. 

\item If $p\mid \gcd(\xi_0, \dots, \xi_n)$ let 
\[
 \a_p:= \min\{  | x_i |_p     \;  | \; \text{ such that } \;   x_i\neq 0 \; \text{ and } \; i=0, \ldots , n.       \}.
 \]
%
Then $f^M$ is semistable over $K_p$ for  $M=\begin{bmatrix}  \frac 1 {p^{r_p}} & 0 \\ 0 & 1 \end{bmatrix}$ and   $r_p= \frac {2{\a_p}} {d q_j}$.
 \end{enumerate}
\end{lem} 
 
\proof
From \cref{semistable-0}, a binary form      $f$ is semistable if and only if there exists some $\xi_j$ such that $\xi_j\neq 0$ in $K_p$.  Hence, the first claim of the theorem.  

Assume $p\mid \xi_i$, for all $i=0, \dots n$. We can further assume that $\wgcd (\xi (f) )=1$ so $f$ is minimal as in \cite[Prop.~6]{b-g-sh}.
Pick $\xi_j$ such that $\xi_j\neq 0$.
Let   $ \xi_j = p^\a \beta$ such that $\gcd (\a, \beta )=1$ and take 
\[
M= \begin{bmatrix}  \frac 1 {p^r} & 0 \\ 0 & 1 \end{bmatrix} \in \GL_2 (\bar K)
\]
for $r= \frac {2\a} {d q_j}$.  Then  $f^M (x, y) = f \left( \frac x {p^r}, y    \right) $ and  from \cref{lem-1} we have 
\[
\xi (f^M) = \left[   \left( \frac 1 {p^r}  \right)^{\frac 1 2 d q_0}  \xi_0,   \dots ,      \b , \dots 
 \left( \frac 1 {p^r}  \right)^{\frac 1 2 d q_n}  \xi_n       \right] 
\]
and $p \nmid \b$. This completes the proof. 
\qed

\subsubsection{Stability over global fields}
Next we want to give a criteria for stability over global fields.  A point $\p=\xi (f)\in {\mathbb{WP}}_\w^n (K)$ is unstable if there is a prime $p\in {\mathcal O}_K$ such that $p \mid \gcd (\xi_0, \dots , \xi_n)$.    Assume there is such a $p \mid \gcd (\xi_0, \dots , \xi_n)$. 

\begin{prop}
Let $ f \in {\mathcal O}_K [x, y]$ and      $\p=\xi(f)=[\xi_0, \dots , \xi_n]  \in {\mathbb{WP}}_\w^n ({\mathcal O}_K)$.  Assume $\p$ is normalized (i.e. multiply  $\p$ by $\frac 1 {\wgcd(\p)}$). 
Then $f^M$ is semistable over $K$,  for $M= \begin{bmatrix} \frac 1 \lambda & 0 \\ 0 & 1 \end{bmatrix}$ and 
\[
\lambda = \prod_{p \mid \gcd (\xi_0, \dots , \xi_n) } p^{r_p}
\]
\end{prop}

\proof
If $\p$ is normalized then $\wgcd \left(\xi_0(f), \dots , \xi_n(f) \right) =1$. Let  
\[
\gcd \left(\xi_0(f), \dots , \xi_n(f) \right) = \prod_{i=1}^s p_i^{a_i},
\]
 where  
$p_i \in {\mathcal O}_K$ are  primes.    Then from the above Lemma, exists $r_i$ such that for 
$M_i = \begin{bmatrix} \frac 1 {\lambda_i} & 0 \\ 0 & 1 \end{bmatrix}$ 
the form $f^{M_i}$ is semistable over $K_{p_i}$.  Let $M= \prod_{i=1}^s M_i$.  Then $f^M$ is semistable for every prime $p_i \mid \gcd (\xi_0(f), \dots , \xi_n(f) )$,  hence it is semistable over $K$.
\qed

\section{Semistability and weighted heights} \label{sec-4}


%

Let $K$ be a number field and ${\mathcal O}_K$ its ring of integers. $M_K$ denotes the set of places of $K$, where $M_K$ are the Archimedean places and $M_\infty$ the non Archimedean places of $K$.  For any norm $\nu \in M_K$, the completion of $K$ at $\nu$ is denoted by $K_\nu$.

The group action $K^\star$ on ${\mathbb A}^{n+1} (K)$ induces a group action of ${\mathcal O}_K$ on ${\mathbb A}^{n+1} (K)$.         
By  ${\mathcal O}rb (\p)$ we denote the ${\mathcal O}_K$-orbit in ${\mathbb A}^{n+1} ({\mathcal O}_K)$ which contains $\p$. 
For any point $\p =[x_0 : \dots : x_n] \in {\mathbb{WP}}_w^n (K)$   we can assume, without loss of generality, that $\p =[x_0 : \dots : x_n] \in {\mathbb{WP}}_w^n ({\mathcal O}_K)$. The  height for weighted projective spaces will be defined in the next section.  

For the rest of this section we assume $K={\mathbb Q}$.   
For the tuple $\x = (x_0, \dots , x_n) \in {\mathbb Z}^{n+1}$ we define   the \textbf{weighted greatest common divisor} with respect to the absolute value $| \, \cdot \, |_v$, denoted by $\wgcd_v (\x)$,   
\[ \wgcd_v (\x) :=  \prod_{  \stackrel {d^{q_i} \mid x_i} {d\in {\mathbb Z}}    } | d |_v \]
as the product of all divisors  $d\in {\mathbb Z}$ such that for all $i=0, \dots , n$, we have $d^i \mid x_i$.  
We will call a point $\p \in {\mathbb{WP}}_{\w}^n ({\mathbb Q})$ \textbf{normalized} if $\wgcd (\p) =1$.

\begin{defi}\label{height}
Let $\w=(q_0, \dots , q_n)$ be a set of weights and ${\mathbb{WP}}^n (K)$ the weighted  projective space  over  a number field $K$.   Let  $\p \in {\mathbb{WP}}^n(K)$ a point such that  $\p=[x_0, \dots , x_n]$. We define the  \textbf{weighted multiplicative height} of $\p$   as  
\begin{equation}\label{def:height}
\wh_K( \p ) := \prod_{v \in M_K} \max   \left\{   \frac{}{}   |x_0|_v^{\frac {n_v} {q_0}} , \dots, |x_n|_v^{\frac {n_v} {q_n}} \right\}
\end{equation}
The \textbf{logarithmic height} of the point $\p$ is defined as follows
\begin{equation}
\wh^\prime_K(\p) := \log \wh_K (\p)=   \sum_{v \in M_K}   \max_{0 \leq j \leq n}\left\{\frac{n_v}{q_j} \cdot  \log  |x_j|_v \right\}.
\end{equation}
\end{defi}
%
For any  $\p\in {\mathbb{WP}}_\w^n (K)$:
\begin{enumerate}[\upshape(i), nolistsep]
\item  $\wh_K(\p)$ is well defined
\item $\wh_K(\p) \geq 1$,
\end{enumerate}
see \cite{b-g-sh}. 
The height of   $\p \in {\mathbb{WP}}^n(\overline {\mathbb Q})$ is called the \textbf{absolute (multiplicative) weighted height} and is the function 
\[
\begin{split}
\awh: {\mathbb{WP}}^n(\bar {\mathbb Q}) & \to [1, \infty)\\
\awh(\p)&=\wh_K(\p)^{1/[K:{\mathbb Q}]},
\end{split}
\]
where  $\p \in {\mathbb{WP}}^n(K)$, for any $K$ which contains ${\mathbb Q} (\awgcd (\p))$.  The \textbf{absolute (logarithmic) weighted height} on ${\mathbb{WP}}^n(\overline {\mathbb Q})$  is the function 
\[
\begin{split}
\awh^\prime: {\mathbb{WP}}^n(\bar {\mathbb Q}) & \to [0, \infty)\\
\awh^\prime (\p)&= \log \, \wh (\p)= \frac 1 {[K:{\mathbb Q}]} \awh_K(\p).
\end{split}
\]

%
\subsection{Invariant height}
Let  $v\in M_K$ be a place and $\xi (f)$ the set of invariants of a degree $d$ binary form $f \in V_d$. 
 We define the norm 
 \[
 {\mid {\xi ( f )} \mid} := \max_{0\leq i \leq n} \left\{  \norm{\xi_i}_v^{\frac 1 {q_i}}   \right\}, 
 \]
 and 
\[
\norm{\xi }_v^t ( f  )  := \frac {\mid {\xi ( f )} \mid ^t} {\max_i \{   \mid {f_i} \mid _v^t \}  }
\]
see \cite{rabina} or \cite{r-savin} for details. 

Let $D$ be the divisor determined by the roots of $f(x, y)$ via the Chow coordinates.    The invariant height of the divisor $D$ as defined in \cite{zhang} is 
\[
h (D)  = \left(     \prod_{v \in M_K}  \inf_{M \in \SL_2 ({\mathbb C})}  \norm{ \xi }_v^t \left(  f^M  \right)     \right)^{\frac  1 {[K:{\mathbb Q}]}}
\]
The \textbf{logarithmic invariant height} of   $D$ is defined as  
\[
\hat h (D) = \frac 1 {[K:{\mathbb Q}]} \sum_{v \in M_K}   \inf_{M \in \SL_2 ({\mathbb C})} \left(  \frac { - \log \norm{\xi }_v^t \left(  f^M  \right) }   t   \right)  
\]
%
 
\begin{rem} In \cite[Thm.4.4.2]{rabina} it is claimed that 
if $f$ is a degree $\deg f=d$  semi-stable binary form then $\hat h(f)\geq 0$. Moreover, if ${\Delta}_f \neq 0$ then 
\[
\hat h (f)  \geq \frac d 2 \log \frac 4 3 
\]
If $f$ is a stable binary form then $\hat h(f)> 0$.
\end{rem}

 In \cite{rabina} it is claimed that for $f(x, y) = x^d - y^d$, for $d\neq p^r$, for any prime $p\in \O_K$ then 
 \[
 \hat h (f) =  \frac d 2  \log 2
 \]


On contrary to the invariant height in \cite{zhang} the weighted moduli height is very easily computed. 
Next we will determine the heights of strictly semistable and semistable points in ${\mathbb{WP}}_\w^n (K)$.

\begin{thm}\label{thm-6}
Let $f \in V_d$ such that   $f \in {\mathcal O}_K [x, y]$. Then the following hold:

\begin{enumerate}[\upshape(i), nolistsep]
\item  If $f$ is semistable then its  weighted moduli height    $\wh (\xi(f)) \geq 1$.
\item  If $f$ is strictly semistable then  $d$ is even and its weighted moduli height $\wh (\xi(f)$ is determined in  \cref{tab-1} for   $d=4, 6, 8, 10$.
\end{enumerate}
\end{thm}

\begin{table}[htp]       
\caption{Strictly semistable points and their weighted heights}
\begin{center}
\begin{tabular}{|c|c|c|c|}
\hline
d     &  $\xi(f)$ &  $\wh (\xi (f))$  &  $\awh^\prime (\xi (f))$ \\
\hline 
4     & $ [1: -2]$  & $2 ^ {1/3}$ &  0.231 \\
\hline 
6        & $[-3 : 3 :- 1 :-3^5] $ & $ 3\sqrt3 $ & 0.549 \\
\hline
8   & $[ 2: 2^2 \cdot 3 : 2^6 : 2^6 : 2^9 : 2^9]$  & $2 \sqrt2$ & 1.039\\
\hline
10    & $ [- 5 , 5^4, -2^2 5^{7}, -2^2 5^4, 5^{8}, 0, -2^3 5^{11}, -2^2 5^{7}, -2^3 5^{15}]$ & $5 \cdot 2^{1/3}  \cdot 5^{1/6}  $   &    2.108  \\
\hline
\end{tabular}
\end{center}
\label{tab-1}
\end{table}

\proof   Part i) is a direct consequence of the definition of the weighted height in \cite{b-g-sh}.  


Let $f$ be a binary quartic which is strictly semistable.     Then 
\[
\xi(f)=  \left[ \xi_0: \xi_1 \right] = [1: -2].
\] 
    Let $f$ be a sextic with a root of multiplicity 3.  
Invariants of $f$ with weights $\w=(2,4,6,10)$ are given by  
\[ 
\left[ \xi_0: \xi_1: \xi_2:  \xi_3 \right] = \left[ -3 : 3 : -1 : -3 ^ 5 \right].
\]
An octavic $f(x, y)$ is strictly semistable  if and only if the basic invariants with weights $\w=(2,3,4, 5, 6, 7)$ take the form
\[
 \xi (f) = [\xi_0: \xi_1: \xi_2:  \xi_3: \xi_4: \xi_5] =
  [   2,  2^2 \cdot 3,  2^6,  2^6,  2^9,  2^9 ]
\]
For the decimic binary form
\[
f(x, y) = x^5 \left(x^5+x^4 y a_4+x^3 y^2 a_3+x^2 y^3 a_2+x \,y^4 a_1+y^5\right) 
\]
the strictly semistable point is given below:
\[
\begin{split}
\xi ( f ) & =  \left[ \xi_0: \xi_1: \xi_2:  \xi_3: \xi_4: \xi_5: \xi_6:  \xi_7: \xi_8 \right]  \\
& = [- 5 , 5^4, -2^2 5^{7}, -2^2 5^4, 5^{8}, 0, -2^3 5^{11}, -2^2 5^{7}, -2^3 5^{15}]
\end{split}
\]
This completes the proof. 

\qed

In the fourth column of \cref{tab-1} we have presented the logarithmic weighted height.  It seems that such logarithmic height increases steadily as $d$ increases.  It would be interesting to determine how fast the logarithmic height increases and how does it compare to the invariant height (which is also a logarithmic height) as defined in \cite{zhang}.    Discussing in detail how the invariant height compares to the weighted moduli height is the focus of \cite{c-sh-2}.



\bibliographystyle{amsplain}

\bibliography{elira-2021}{}

\end{document}